\documentclass[11pt, a4paper, reqno, numberwithinsect,thm-restate,autoref]{amsart}
\usepackage{amssymb}
\usepackage[colorlinks, bookmarks=true]{hyperref}
\usepackage{amsthm, thmtools}
\usepackage{enumerate, mathrsfs,color,amssymb,mathabx}
\usepackage{verbatim}
\usepackage[T1]{fontenc}
\usepackage{mathptmx}

\newtheorem{theorem}{Theorem}[section]
\newtheorem{proposition}[theorem]{Proposition}
\newtheorem{lemma}[theorem]{Lemma}

\theoremstyle{definition}

\newtheorem{corollary}[theorem]{Corollary}

\theoremstyle{remark}
\newtheorem{remark}[theorem]{Remark}

\numberwithin{equation}{section}
\def\sqr#1#2{{\,\vcenter{\vbox{\hrule height.#2pt\hbox{\vrule width.#2pt
height#1pt \kern#1pt\vrule width.#2pt}\hrule height.#2pt}}\,}}

\usepackage{cleveref}

\begin{document}
\title[Additive mappings preserving orthogonality]{An algebraic characterization of linearity for additive maps preserving orthogonality}

\author[L. Li]{Lei Li}
\address[L. Li]{School of Mathematical Sciences and LPMC, Nankai University, 300071 Tianjin, China.}
\email{leilee@nankai.edu.cn}

\author[S. Liu]{Siyu Liu}
\address[S. Liu]{School of Mathematical Sciences and LPMC, Nankai University, 300071 Tianjin, China.}
\email{760659676@qq.com}

\author[A.M. Peralta]{Antonio M. Peralta}
\address[A.M. Peralta]{Instituto de Matem{\'a}ticas de la Universidad de Granada (IMAG), Departamento de An{\'a}lisis Matem{\'a}tico, Facultad de
	Ciencias, Universidad de Granada, 18071 Granada, Spain.}
\email{aperalta@ugr.es}

\subjclass[2010]{Primary 46B20, 46C05, 39B55, 46C99}
\keywords{Birkhoff orthogonality; Euclidean orthogonality; orthogonality preserving additive mappings; complex inner product spaces} 

\begin{abstract} We study when an additive mapping preserving orthogonality between two complex inner product spaces is automatically complex-linear or conjugate-linear. Concretely, let $H$ and $K$ be complex inner product spaces with dim$(H)\geq 2$, and let $A: H\to K$ be an  additive map preserving orthogonality. We obtain that $A$ is zero or a positive scalar multiple of a real-linear isometry from $H$ into $K$. We further prove that the following statements are equivalent:
	\begin{enumerate}[$(a)$]
		\item $A$ is complex-linear or conjugate-linear.
		\item For every $z\in H$ we have $A(i z) \in \{\pm i A(z)\}$.
		\item There exists a non-zero point $z\in H$ such that $A(i z)  \in \{\pm i A(z)\}$.
		\item There exists a non-zero point $z\in H$ such that $i A(z) \in A(H)$.
	\end{enumerate} 
	\noindent The mapping $A$ neither is complex-linear nor conjugate-linear if, and only if, there exists a non-zero $x\in H$ such that $i A(x)\notin A(H)$ (equivalently, for every non-zero $x\in H$, $i A(x)\notin A(H)$). \smallskip
	
\noindent Among the consequences we show that, under the hypothesis above, the mapping $A$ is automatically complex-linear or conjugate-linear if $A$ has dense range, or if $H$ and $K$ are finite dimensional with dim$(K)< 2\hbox{dim}(H)$.
\end{abstract}

\maketitle

\section{A brief state-of-the-art}

In a recent contribution we have established that every additive mapping with dense image between two complex inner product spaces enjoying the additional property that it preserves orthogonality must be complex-linear or conjugate linear, and in fact a positive scalar multiple of a complex-linear or a conjugate-linear isometry between the spaces. It is natural to ask whether the norm-density of the image is a necessary condition for our conclusion. As we shall see below there are examples of real-linear maps which are orthogonality preserving, but are not complex-linear nor conjugate-linear. Our main goal here is to find a complete characterization of those additive orthogonality preserving maps between complex inner product spaces which are complex-linear or conjugate-linear. Our problem is actually related to other results on linear or additive orthogonality preservers.\smallskip

The so-called Koldobsky-Blanco-Turnšek theorem is a fascinating result asserting that every linear map between normed spaces which preserves orthogonality in the Birkhoff-James sense is necessarily a scalar multiple of an isometry (see \cite{Kol93, BlancoTurnsek}). Needless to say how influencing has been this result after its appearance. Actually, thirteen years separate the results in the setting of real spaces from the conclusion in the complex setting. Let us recall that a vectors $x,y$ in a normed space $X$ over $\mathbb{F}= \mathbb{R}$ or $\mathbb{C}$ are said to be \textit{orthogonal in the Birkhoff-James sense} or \emph{Birkhoff-James orthogonal} ($x\perp_B y$ in short), if for every scalar $\alpha \in\mathbb{K}$ we have $\|x+\alpha y\|\ge \|x\|$. When $X$ is an inner product space (with inner product $\langle\cdot, \cdot\rangle$), $x$ and $y$ are Birkhoff-James orthogonal if, and only if, they are orthogonal in the Euclidean sense, that is, $x\perp_B y \Leftrightarrow x\perp_2 y \Leftrightarrow \langle x|y\rangle =0$. In this note we shall mainly deal with complex inner product spaces, so our notion of orthogonality is the Euclidean orthogonality. A (non-necessarily linear nor continuous) mapping $T: X\to Y$ between normed spaces $X$ and $Y$ is called \textit{preserves orthogonality} if $x\perp_B y$ implies $T(x)\perp_B T(y)$ for any $x,y\in X$. The the equivalence $x\perp_B y \leftrightarrow T(x)\perp_B T(y)$ holds we say $T$ \emph{preserves orthogonality in both directions}. We note that a result by Chmieli\'{n}ski from 2005 \cite{Chmielinski2005} advanced the conclusion in the Blanco-T theorem in the setting of real or complex inner product spaces.\smallskip

Recent studies consider merely additive mappings preserving orthogonality. A remarkable result by W\'{o}jcik in this line reads as follows:

\begin{theorem}$($\cite[Theorem 3.1]{Wojcik2019}$)$
Let $X, Y$ be two real normed spaces. For a non-vanishing additive mapping $A: X\to Y$ the following conditions are equivalent:
\begin{enumerate}[$(a)$]
\item  $A$ preserves $B$-orthogonality, i.e., $\forall\,x,y\in X, x\perp_B y \Rightarrow A(x)\perp_B A(y)$;
\item $A$ is real-linear and there exists $\gamma>0$ such that $\|A(x)\|=\gamma\| x\|$ for any $x\in X$.
\end{enumerate}
\end{theorem}

Motivated by the study of additive preservers of truncations between Cartan factors and atomic JBW$^*$-triples (cf. \cite[\S 3]{LiLiuPe2024addpreserverstruncations}), we were attracted by the natural question of determining whether complex-linearity or conjugate-linearity is an automatic property of additive maps preserving orthogonality between two complex inner product spaces (see \cite{LLP2024AOP}). Surprisingly, a topological condition on an additive preserver of orthogonality between complex inner product spaces like having norm-dense range, suffices to guarantee that it is linear. 

\begin{theorem}\label{t: LiLiuPe2025}$($\cite[Theorem 2.1]{LLP2024AOP}$)$ Let $A: H\to K$ be a non-zero mapping between two complex inner product spaces with $\dim (H)\ge 2$. Suppose that $A$ has dense image. Then the following statements are equivalent:
\begin{itemize}
\item[(a)] $A$ is $($complex$)$ linear or conjugate-linear mapping and there exists $\gamma>0$ such that $\|A(x)\|=\gamma\| x\|$, for all $x\in X$, that is, $A$ is a positive scalar multiple of a linear or a conjugate-linear isometry;
\item[(b)] There exists $\gamma_1>0$ such that one of the next properties holds for all $x,y\in H$:
\begin{itemize}
\item[(b.1)] $\langle A(x)|A(y)\rangle =\gamma_1 \langle x|y\rangle$,
\item[(b.2)] $\langle A(x)|A(y)\rangle =\gamma_1 \langle y|x\rangle$;
\end{itemize}
\item[(c)] $A$ is linear or conjugate-linear and preserves orthogonality in both directions;
\item[(d)] $A$ is linear or conjugate-linear and preserves orthogonality;
\item[(e)] $A$ is additive and preserves orthogonality in both directions;
\item[(f)] $A$ is additive and preserves orthogonality.
\end{itemize}
\end{theorem}

In particular if $H$ and $K$ are complex inner product spaces with dim$(H)\geq 2$, every non-zero surjective additive mapping $A: H\to K$ preserving orthogonality is a positive scalar multiple of a complex-linear or a conjugate-linear isometry (and it is, in particular, continuous). It seems quite natural to ask whether the hypothesis $A$ being surjective or $A$ having norm-dense range, can be relaxed. As we shall see below, there are examples of additive maps preserving orthogonality which are not complex-linear nor conjugate-linear. \smallskip

The purpose of this note is to find necessary and sufficient conditions to guarantee that an orthogonality preserving additive map $A: H\to K,$ where $H$ and $K$ are complex inner product spaces with dim$(H)\geq 2$, is complex-linear or conjugate-linear. We prove, in \Cref{t: characterization of comple linear OP}, that the following statements are equivalent:
\begin{enumerate}[$(a)$]
	\item $A$ is complex-linear or conjugate-linear.
	\item For every $z\in H$ we have $A(i z)  \in \{\pm i A(z)\}$.
	\item There exists a non-zero point $z\in H$ such that $A(i z)  \in \{\pm i A(z)\}$.
	\item There exists a non-zero point $z\in H$ such that $i A(z) \in A(H)$.
\end{enumerate} 
\noindent Moreover, the mapping $A$ neither is complex-linear nor conjugate-linear if, and only if, there exists a non-zero $x\in H$ such that $i A(x)\notin A(H)$ (equivalently, for every non-zero $x\in H$, $i A(x)\notin A(H)$).\smallskip

\Cref{c: if A has norm-dense range every point in the domain must be of mixed type} shows that every non-zero additive mapping $A: H\to K$ preserving Euclidean orthogonality and having norm-dense range satisfies statement $(b)$ above, and thus it must be complex-linear or conjugate-linear. We rediscover in this way \Cref{t: LiLiuPe2025} above. Another consequence of our main result establishes that in case that $H$ and $K$ are finite dimensional complex Hilbert spaces with $2\leq$dim$(K)< 2\hbox{dim}(H)$, every additive mapping $A: H\to K$ preserving orthogonality must be complex-linear or conjugate-linear (see \Cref{c automatic complex-linearity or conjugate linearity for finite ranges}).\smallskip

By relying on W\'{o}jcik's theorem above, we show how the problem considered in our main goal can be reduced to the case of real-linear isometries between complex inner product spaces, since under the above hypotheses, for every non-zero additive mapping $A: H\to K$ preserving Euclidean orthogonality, there exists a positive $\gamma\in\mathbb{R}$ and a  real-linear isometry $T: H\to K$ such that $A = \gamma T$ (see \Cref{p real-linearity and distance preservation}).\smallskip

It is worth to note, thought widely known, that there are examples of surjective real-linear isometries between complex Hilbert spaces which are not complex-linear or conjugate linear and do not preserve orthogonality, for example, $T: \ell_2^2(\mathbb{C}) \to  \ell_2^2(\mathbb{C})$, $T(\alpha,\beta) = (\overline{\alpha},\beta)$.\smallskip

We conclude this note with a couple of results of independent interest; in the first one, we show that every additive orthogonality preserving mapping between two $($real or complex$)$ inner product spaces admits a unique extension to an additive orthogonality preserving map between the corresponding completions of the spaces (see \Cref{l op lift to the completion}). This is finally employed to show that if $A: H\to K$ is an additive orthogonality preserving mapping between two complex inner product spaces, and $\tilde{K}$ denotes the completion of ${K}$, we can always find a bounded real-linear mapping $R: \tilde{K}\to \tilde{K}$ satisfying that $R T: {H}\to \tilde{K}$ is a complex-linear positive scalar multiple of a complex-linear isometry preserving orthogonality (see \Cref{c: final}). 

\section{Real-linearity and automatic continuity}

Our main goal will consist in characterize additive maps preserving orthogonality between complex inner product spaces without extra assumptions on their images. Related these weaker assumptions, in \cite[Proposition 2.2]{LLP2024AOP} we prove that for every couple of normed spaces $X$ and $Y$ such that dim$(X)\geq 2$ and $X$ admits a conjugation (i.e. a conjugate-linear isometry of period-$2$) every additive mapping $A:X\to Y$ preserving Birkhoff orthogonality is real-linear. However, it is not clear why every complex inner product space admits a conjugation. \smallskip

Our first goal is to prove the automatic real-linearity and continuity of all additive orthogonality preserving maps between complex inner product spaces. 

\begin{proposition}\label{p real-linearity and distance preservation} Let $H$ and $K$ be complex inner product spaces such that $\dim (H)\ge 2$. Suppose that $A: H\to K$ is a non-zero additive mapping preserving Euclidean orthogonality. Then there exists a positive $\gamma\in\mathbb{R}$ and a real-linear isometry $T: H\to K$ such that $A = \gamma T$. 
\end{proposition}

\begin{proof} As we have commented above, if $H$ and $K$ are Hilbert spaces, the real-linearity of $A$ follows from \cite[Proposition 2.2]{LLP2024AOP}. In the case of inner product spaces we exploit an idea in the proof of \cite[Theorem 2.1]{LLP2024AOP}.  Since $A\neq 0$, there exists $x_0\in H$ such that $A(x_0)\neq 0$. For each $x_1\in H\backslash\{ 0\}$ satisfying $x_1\perp_2 x_0$ in $H$ (the existence of such an $x_1$ is guaranteed because $\dim(H)\ge 2$). Let $H_{\mathbb{R}}(x_0, x_1)=\mathbb{R} x_0\oplus \mathbb{R} x_1$ denote the (2-dimensional) real subspace of $H$ generated by the vectors $\{x_0,x_1\}.$ \smallskip

\noindent\emph{Claim 1} For each $x_a\in H$, $x_1\in  H\backslash\{ 0\}$ with $x_1\perp_2 x_a$ and $A(x_a)\neq 0$, the restricted mapping $A|_{H_{\mathbb{R}}(x_a, x_1)} : H_{\mathbb{R}}(x_a, x_1)\to K$ is a positive scalar multiple of a real-linear isometry. Furthermore, taking $\gamma (x_a) = \frac{\left\|A|_{H_{\mathbb{R}(x_a,x_1)}}(x_a)\right\|}{\|x_a\|} \in \mathbb{R}^+$, the mapping $\gamma (x_a)^{-1} A|_{H_{\mathbb{R}}(x_a, x_1)}$ is a real-linear isometry. In particular $A(x_1)\neq 0$.  \smallskip
	
To prove the claim we first observe that $H_{\mathbb{R}}(x_a, x_1)$ is a real Hilbert space for the inner product of $H$. Actually, for any $x=\lambda_0 x_a+\lambda_1 x_1, y=\mu_0 x_a+\mu_1 x_1$ in $H_{\mathbb{R}}(x_a,x_1),$  where $\lambda_0, \lambda_1,\mu_0,\mu_1\in\mathbb{R}$, since
\[\langle x |y\rangle =\lambda_0\mu_0\|x_a\|^2+\lambda_1\mu_1\|x_1\|^2\in\mathbb{R},\]
we can see that $\langle x|y\rangle=\Re\hbox{e}\langle x|y\rangle$ for all $x,y\in H_{\mathbb{R}}(x_a,x_1)$, and thus elements in $H_{\mathbb{R}}(x_a,x_1)$ are orthogonal in $(H, \langle \cdot |\cdot \rangle)$ if, and only if, they are orthogonal in the real Hilbert space $(H_{\mathbb{R}}(x_a,x_1),\Re\hbox{e} \langle \cdot |\cdot\rangle)$. Let $K_{\mathbb{R}}$ denote the real inner product space underlying $K$ with inner product given by $(x |y):=\Re\hbox{e}\langle x |  y \rangle$ ($x,y\in K$).  Since $x_a\in H_{\mathbb{R}}(x_a,x_1)$ and $A(x_a)\neq 0$, we can therefore conclude that the mapping $$A|_{H_{\mathbb{R}}(x_a,x_1)}:H_{\mathbb{R}}(x_a,x_1)\to K_{\mathbb{R}}$$ also is a non-zero additive mapping preserving Euclidean orthogonality between two real inner product spaces. According to W\'{o}jcik's theorem (see \cite[Theorem 3.1]{Wojcik2019}), there exists $\gamma(x_a,x_1)\in\mathbb{R}^+$, depending on $x_a$ and $x_1$, such that $A|_{H_{\mathbb{R}}(x_a,x_1)}$ is real-linear and $\|A(x)\|=\gamma(x_a,x_1)\| x\|$ for all $x$ in $H_{\mathbb{R}}(x_a,x_1)$. Since 
$\|A|_{H_{\mathbb{R}}(x_a,x_1)}(x_a)\|=\gamma(x_a,x_1)\|x_a\|$, it follows that  $\gamma(x_a,x_1)=\frac{\|A|_{H_{\mathbb{R}(x_a,x_1)}}(x_a)\|}{\|x_a\|}$, which does not depend on $x_1$ and will be denoted by $\gamma(x_a)$. It also follows that $\| A(x)\| = \gamma(x_a) \|x\|$ for all $x\in H_{\mathbb{R}}(x_a,x_1)$, and thus $A(x_1) \neq 0$.\smallskip

\noindent\emph{Claim 2} For each $x_1\in  H\backslash\{ 0\}$ with $x_1\perp_2 x_0$ and $\alpha\in \mathbb{C}\setminus \{0\}$ we have \begin{equation}\label{gamma}
	\frac{\|A(\alpha x_1)\|}{\| \alpha x_1\|}=\frac{\|A(\alpha x_0)\|}{\| \alpha x_0\|}=\frac{\|A(x_1)\|}{\|x_1\|}=\frac{\|A(x_0)\|}{\|x_0\|}.
\end{equation} Observe that $\alpha x_1\in H\setminus \{0\}$ and $\langle \alpha x_1 |x_0\rangle =0$. So, by applying the conclusion in \emph{Claim 1} to the subspaces $H_{\mathbb{R}}(x_0,\alpha x_1)$ and $H_{\mathbb{R}}(x_0, x_1)$ we get $ \frac{\|A(\alpha x_1)\|}{\| \alpha x_1\|}=\frac{\|A(x_0)\|}{\|x_0\|} = \frac{\|A(x_1)\|}{\| x_1\|}\neq 0$. If in \emph{Claim 1} we replace $x_0$ and $x_1$ with $x_1$ and $\alpha x_0$, respectively, we get $\frac{\|A(\alpha x_0)\|}{\| \alpha x_0\|}=\frac{\|A(x_1)\|}{\|x_1\|}$.\smallskip

\noindent\emph{Claim 3} For each $z\in  H\backslash\{ 0\}$ we have \begin{equation}\label{gamma} \frac{\|A(z)\|}{\|z\|}=\frac{\|A(x_0)\|}{\|x_0\|},
\end{equation} that is, $\frac{\|x_0\|}{\|A(x_0)\|} A$ is an isometry.\smallskip

Since $\dim H\ge 2$, there exists $x_1\neq 0$ satisfying $x_1\perp_2 x_0$ such that $z=\alpha x_0+\beta x_1$, with $\alpha, \beta \in \mathbb{C}$.\smallskip

If $\alpha=0$, the element $z=\beta x_1$ is non-zero and it follows from \emph{Claim 2} (see \eqref{gamma}) that 
$$\frac{\|A(z)\|}{\| z\|}= \frac{\|A(\beta x_1)\|}{\| \beta x_1\|}=\frac{\|A(x_0)\|}{\|x_0\|}.$$

If $\beta=0$, the element $z=\alpha x_0$ is non-zero, and we get from \emph{Claim 2} (see \eqref{gamma}) that $$\frac{\|A(z)\|}{\| z\|}= \frac{\|A(\alpha x_0)\|}{\| \alpha x_0\|}=\frac{\|A(x_0)\|}{\|x_0\|}.$$

Assume next that $\alpha, \beta \neq 0$. Since $\alpha x_0\perp \beta x_1$, by applying \emph{Claim 1} to the subspace $H_{\mathbb{R}}(\alpha x_0,\beta x_1)$ and \emph{Claim 1} we derive that
$$\frac{\|A(z)\|}{\| z\|}= \frac{\|A(\alpha x_0+\beta x_1)\|}{\| \alpha x_0+\beta x_1\|}=\frac{\|A(\alpha x_0)\|}{\|\alpha x_0\|},$$ which concludes the proof of \emph{Claim 3}. \smallskip

Finally, if we take $\gamma= \frac{\|A(x_0)\|}{\|x_0\|}\in \mathbb{R}^+$, it follows from \emph{Claim 3} that the mapping $\gamma^{-1} A: H\to K$ is an additive isometry, and thus it is real-linear, as desired. 
\end{proof}

Our next goal is to characterize when a non-zero additive mapping $A$ preserving orthogonality between complex inner product spaces $H$ and $K$ is actually complex-linear or conjugate linear. We shall see below that this is not always the case. By \Cref{p real-linearity and distance preservation} we can find $\gamma>0$ and a real-linear isometry $T:H\to K$ satisfying $A = \gamma T$. Furthermore, $T$ preserves orthogonality because $A$ does. So, we can focus on real-linear isometries preserving orthogonality.

\begin{lemma}\label{l: image of a norm-one vector x0} Let $H$ and $K$ be complex inner product spaces with dim$(H)\geq 2$, and let $T: H\to K$ be a real-linear isometry preserving orthogonality. Then for each norm-one vector $x_0\in H$, there exist a unique $\alpha(x_0) \in i\mathbb{R}$ and a unique $\eta(x_0) \in K$ $($both depending on $x_0)$ satisfying:
\begin{enumerate}[$(1)$]
\item $\langle \eta(x_0) | T(x_0)\rangle =0$, that is $T(x_0)\perp \eta(x_0)$,
\item $|\alpha(x_0)|^2+\|\eta(x_0)\|^2=\| x_0\|^2=1$,
\item $T(i x_0)=\alpha(x_0) T(x_0)+\eta(x_0)$.
\end{enumerate} 
\end{lemma}

\begin{proof} 
Observe that $\|T(x_0)\|=\|x_0\|=1$ and $\|T(i x_0)\|= \|i x_0\|=1$ (that is, $\{ T(x_0), T(i x_0)\}$ is contained in the unit sphere of $K$). By defining $\eta(x_0) := T (i x_0)- \langle T(i x_0) | T(x_0)\rangle T(x_0)\in K$ and $\alpha(x_0) := \langle T(i x_0) | T(x_0)\rangle\in \mathbb{C}$,  the identity $T(i x_0)= \alpha(x_0) T(x_0)+\eta(x_0)$ clearly holds, and by construction $\langle \eta(x_0) | T(x_0)\rangle =0$. Moreover, by definition we also have 
\[1=\|T(ix_0)\|^2=|\alpha(x_0) |^2\|T(x_0)\|^2+\|\eta(x_0)\|^2=|\alpha(x_0)|^2+\|\eta(x_0)\|^2.\]

Pick $c,s\in\mathbb{R}$ such that $c^2+s^2=1$. By considering $c x_0 +s i x_0$ as an element in the unit sphere of $H$, and having in mind that $T$ is a real-linear isometry, we deduce that $T(c x_0 + s i x_0)=c T(x_0)+ s T(i x_0) = c T(x_0) + s \alpha(x_0) T(x_0) + s \eta(x_0)$ must have norm-one. Therefore 
\begin{equation}
\begin{split}
1&=\|T( c x_0 + s i x_0 )\|^2 =\|(c+s \alpha(x_0)) T(x_0) + s\eta(x_0) \|^2  \\ & =|c+s \alpha(x_0) |^2 + s^2 \|\eta(x_0)\|^2 
=c^2 + s^2 |\alpha(x_0)|^2 + 2 c s \Re\hbox{e}(\alpha(x_0)) + s^2 \|\eta(x_0)\|^2 \\
&=c^2 + s^2 + 2 c s \Re\hbox{e} (\alpha(x_0)) =1 + 2 c s \Re\hbox{e} (\alpha(x_0)),
\end{split}
\end{equation} for all $c,s\in \mathbb{R}$ as above, which finally implies that $\Re\hbox{e}(\alpha(x_0))=0,$ and $\alpha(x_0)\in i\mathbb{R}$.
\end{proof}

Under the hypotheses of \Cref{l: image of a norm-one vector x0}, we can actually deduce that for each non-zero $z\in H$ there exist a unique $s(z)\in \mathbb{R}$ (with $-1\leq s(z)\leq 1$) and a unique $\eta (z)\in K$ such that $T(z) \perp \eta (z)$ and  
\[T( i z )=i s(z)  T(z) + \eta (z).\] In particular, $\|\eta (z)\|^2 = (1-s(z)^2) \|z\|^2$. We shall say that $z$ is a \emph{point of pure complex type for $T$} if $\eta (z) =0$, otherwise we say that $z$ is a \emph{point of mixed complex type for $T$}. Note that if $z$ is a point of pure complex type for $T$ we have $s(z) \in \{\pm 1\}$.\smallskip

In the next results we show that a real-linear mapping between two inner product spaces cannot admit two orthogonal directions of pure and mixed complex type respectively. 

\begin{lemma}\label{l: T does not mix different types at orthogonal at orhtogonal directions} Let $H$ and $K$ be complex inner product spaces with dim$(H)\geq 2$, and let $T: H\to K$ be a real-linear isometry preserving orthogonality. Suppose that there exists a norm-one element $x_0\in H$ which is of pure complex type for $T$, that is, $T(i x_0)=i s(x_0) T(x_0)$, with $s(x_0)\in\{\pm 1\}$. Then every norm-one vector $x_1\in H$ with $x_1\perp x_0$ is of pure complex type and $s(x_1) = s(x_0)$, that is, $T(i x_1)=i s(x_0) T(x_1)$.
\end{lemma}

\begin{proof} Arguing by contradiction, we assume that there exists a norm-one element $x_1\in S_H$ such that $x_1\perp x_0$ and $T(i x_1)=i s(x_1) T(x_1) + \eta(x_1)$ with $\eta(x_1)\perp T(x_1)$ and $\eta(x_1)\neq 0$. Since $x_0\perp x_1$, we also have $x_0\perp ix_1$ , and  since $T$ preserves orthogonality, it follows that $\langle T(x_0)|T(x_1)\rangle=0 = \langle T(x_0)|T(i x_1)\rangle.$ We therefore have:
$$ 0=\langle T(x_0) | T(ix_1)\rangle =- i s(x_1) \langle T(x_0)|T(x_1)\rangle + \langle T(x_0)|\eta (x_1)\rangle =\langle T(x_0)|\eta (x_1)\rangle,$$ which proves that $T(x_0)\perp\eta(x_1)$.\smallskip

On the other hand, since the vectors $ix_1+x_0$ and $x_1+ix_0$ are orthogonal in $H$, the hypotheses on $T$ assure that  $T(ix_1+x_0)\perp T(x_1+ix_0)\,\text{ in }\, K,$ which implies:
$$\begin{aligned}
	0&=\langle T(ix_1+x_0)|T(x_1+ix_0)\rangle =\langle T(ix_1)+T(x_0)| T(x_1)+T(ix_0)\rangle \\
	&=\langle i s(x_1) T(x_1) + \eta(x_1) + T(x_0) | T(x_1)+ i s(x_0) T(x_0)\rangle = i s(x_1) - i s(x_0) = i (s(x_1) - s(x_0)).
\end{aligned} $$
Thus, we have $s(x_1)=s(x_0)\in\{\pm 1\}$, and thus $\eta (x_1) = 0$ (since $s(x_1)^2 + \|\eta(x_1)\|^2 =\|T(x_1)\|^2 =\|x_1\|^2=1$), which contradicts the assumption on $x_1$. Observe that the above arguments also show that $s(x_1) = s(x_0)$ for all $x_1$ satisfying the hypotheses of the lemma. 
\end{proof}

Next we see how the existence of a point of mixed complex type in the domain forces that every orthogonal vector to that point is of mixed complex type. 

\begin{proposition}\label{p: mixed complex type point} Let $H$ and $K$ be complex inner product spaces with dim$(H)\geq 2$, and let $T: H\to K$ be a real-linear isometry preserving orthogonality. Suppose that there exists a norm-one element $x_0\in H$ which is of mixed complex type for $T$, that is, $T(i x_0)= i s(x_0) T(x_0) + \eta (x_0)$, with $\eta(x_0)\neq 0$. Then every norm-one element $x_1$ in $H$ which is orthogonal to $x_0$ is of mixed complex type for $T$. Moreover, under these conditions, $s(x_0 ) = s(x_1)$, $\| \eta(x_0) \| = \| \eta (x_1)\|\neq 0,$ and the sets $\{T(x_0), \eta(x_0)\}$ and  $\{T(x_1), \eta(x_1)\}$ are orthogonal in $K$.  
\end{proposition}

\begin{proof} \Cref{l: T does not mix different types at orthogonal at orhtogonal directions} assures that $x_1$ must be of mixed complex type too (that is, $\eta(x_1)\neq 0$). Since $x_0\perp x_1$ and $T$ preserves orthogonality, we derive that $T(x_0), T(i x_0)\perp T(x_1), T(i x_1)$. It follows that 
$$ 0=\langle T(i x_0) | T(x_1)\rangle =i s(x_0) \langle T(x_0) | T(x_1) \rangle + \langle \eta(x_0) | T(x_1)\rangle = \langle \eta(x_0) | T(x_1)\rangle,$$
which implies that $\eta(x_0) \perp T(x_1)$. We can similarly conclude that $\eta(x_1)\perp T(x_0)$.\smallskip

By combining all the previous conclusions we get $$\begin{aligned}
 0 = \langle T(i x_1)|T(i x_0)\rangle =& s(x_0) s(x_1) \langle T(x_1) | T(x_0) \rangle + i s(x_1) \langle T(x_1) |\eta(x_0) \rangle \\
 &-i s(x_0) \langle \eta(x_1) | T(x_0) \rangle + \langle \eta(x_0) |\eta(x_1) \rangle\\
 & = \langle \eta(x_0) |\eta(x_1) \rangle,
\end{aligned}$$ witnessing that $\eta(x_0)\perp \eta(x_1)$, and thus \begin{equation}\label{eq T and eta are orthogonal} \{T(x_0),\eta(x_0)\}\perp \{T(x_1),\eta(x_1)\}.
\end{equation}

Having in mind that the vectors $i x_0 + x_1$ and $x_0 + i x_1$ are orthogonal in $H$, their images $T( i x_0 + x_1 )$ and $T( x_0 + i x_1 )$ are orthogonal in $K$, which combined with \eqref{eq T and eta are orthogonal} leads to $$\begin{aligned}
0&=\langle T(ix_0+x_1)|T(x_0+ix_1)\rangle =\langle T(ix_0) + T(x_1) | T(x_0) + T(ix_1)\rangle\\
&=\langle i s(x_0) T(x_0) + \eta(x_0) + T(x_1) | T(x_0) + i s(x_1) T(x_1) + \eta(x_1)\rangle =i s(x_0) - i s(x_1). 
\end{aligned}
$$ We have therefore shown that $s(x_0)=s(x_1)$. Finally, the equality $s(x_0)^2 + \|\eta(x_0)\|^2 =\|x_0\|^2 = 1= \|x_1\|^2 = s(x_1)^2 +\| \eta(x_1)\|^2$, gives $\| \eta(x_0) \| = \| \eta (x_1)\|$.
\end{proof}

We can now state the key ingredient for our main theorem in this section. Namely, the existence of a point of pure complex type in the domain implies that all the other points are of pure complex type. 

\begin{proposition}\label{p the existence of a pure complex point implies that all are of pure complex type} Let $H$ and $K$ be complex inner product spaces with dim$(H)\geq 2$, and let $T: H\to K$ be a real-linear isometry preserving orthogonality. Suppose that there exists a norm-one element $x_0\in H$ which is of pure complex type for $T$, that is, $T(i x_0)=i s(x_0) T(x_0)$, with $s(x_0)\in\{\pm 1\}$. Then every non-zero vector $x$ in $H$ is of pure complex type too and $s(x_0)=s(x)$, that is, $T(i x ) =i s(x_0) T(x),$ for all $x\in H$.  
\end{proposition}

\begin{proof} If $x\perp x_0$, the conclusion follows from \Cref{l: T does not mix different types at orthogonal at orhtogonal directions}.\smallskip

We assume next that $x\in H$ is non-zero and $\langle x | x_0\rangle \neq 0$. Since dim$(H)\geq 2$, there exist $\alpha, \beta\in \mathbb{C}$ and a norm-one $x_1\in H$ satisfying $x_1\perp x_0$ and $x=\alpha x_0+\beta x_1$. \smallskip

The mentioned \Cref{l: T does not mix different types at orthogonal at orhtogonal directions} asserts that 
\[T(ix_0)=i s(x_0) T(x_0),\quad\text{and}\quad T(i x_1)= i s(x_0) T(x_1).\]
By applying that $T$ is real-linear and the above properties we deduce that 
$$
\begin{aligned}
T(ix)&=T\left(i\Re\hbox{e}(\alpha) x_0\right)+T\left(-\Im\hbox{m} (\alpha) x_0\right)+T\left(i\Re\hbox{e}(\beta) x_1\right)+T\left(-\Im\hbox{m} (\beta) x_1\right)\\
&=\Re\hbox{e} (\alpha) T(ix_0)-\Im\hbox{m} (\alpha) T(x_0) +\Re\hbox{e} (\beta)T(ix_1)-\Im\hbox{m}(\beta) T(x_1)\\
&=\Re\hbox{e} (\alpha) i s(x_0) T(x_0)+\Im\hbox{m}(\alpha) i^2 s(x_0)^2 T(x_0)+\Re\hbox{e} (\beta) i s(x_0) T(x_1) + \Im\hbox{m}(\beta) i^2 s(x_0)^2 T(x_1)\\
&=i s_0(x) \left(\Re\hbox{e} (\alpha) T(x_0) + \Im\hbox{m} (\alpha) T(i x_0) \right)+i s(x_0) \left(\Re\hbox{e} (\beta) T(x_1)+\Im\hbox{m} (\beta) T(i x_1) \right)\\
&=is(x_0) \left(T((\Re\hbox{e} (\alpha)+i\Im\hbox{m} (\alpha))x_0)+T((\Re\hbox{e} (\beta)+i\Im\hbox{m} (\beta))x_1)\right)\\
&=is(x_0) T(\alpha x_0+\beta x_1) =i s(x_0) T(x).
\end{aligned}
$$
\end{proof}

We can now characterize those real-linear isometries preserving orthogonality between complex inner product spaces which are complex-linear or conjugate-linear.

\begin{theorem}\label{t characterization of complex-linear isometries OP} Let $H$ and $K$ be complex inner product spaces with dim$(H)\geq 2$, and let $T: H\to K$ be a real-linear isometry preserving orthogonality. Then the following statements are equivalent:
\begin{enumerate}[$(a)$]
\item $T$ is complex-linear or conjugate-linear.
\item Every non-zero point in $H$ is of pure complex type for $T$.
\item There exists a non-zero point in $H$ which is of pure complex type for $T$.
\item There exists a non-zero point $x_0\in H$ such that $i T(x_0) \in T(H)$.
\end{enumerate}
\end{theorem}

\begin{proof} The implications $(a)\Rightarrow (b) \Rightarrow (c) \Rightarrow (d)$ are clear. We shall prove that $(d) \Rightarrow (a)$. Let us assume the existence of $x_0\in H\setminus \{0\}$ such that $i T(x_0)\in T(H)$. Therefore, there exists a non-zero $x_1\in H$ satisfying $T(x_1) = i T(x_0)$. Clearly, $\|x_0\| = \|T(x_0)\| = \|T(x_1)\| = \|x_1\|$. We can therefore assume that $\|x_1\| = \|x_0\| =1$.\smallskip 
	
\Cref{l: image of a norm-one vector x0} assures that $x_1$ is of pure complex type or of mixed complex type for $T$. We shall proved that the latter is impossible. Suppose, on the contrary, that $x_1$ is a point of mixed complex type for $T$, that is, $T(i x_1) = i s(x_1) T(x_1) + \eta (x_1)$ with $\eta (x_1) \neq 0$ and $\eta(x_1)\perp T(x_1)$. \Cref{p the existence of a pure complex point implies that all are of pure complex type} assures that $i x_1,$ and $x_0$ are of mixed complex type for $T$ with $s(i x_1) = s(x_1) = s(x_0) = s(i x_0)$ and $\|\eta(x_0)\| = \|\eta(i x_0)\|= \|\eta(x_1)\| = \|\eta(i x_1)\|$. By assumptions we have 
$$\begin{aligned}
i T(x_0) &=  T(x_1) = - T(i^2 x_1) = - \Big( i s(i x_1) T(i x_1) + \eta(i x_1) \Big) \\
&= - \Big( i s(x_0) \Big(i s(x_1) T(x_1) + \eta(x_1) \Big) + \eta(i x_1) \Big) \\
&=  s(x_0)^2 T(x_1) -i s(x_0) \eta(x_1) -  \eta(i x_1) = T(x_1) -i  s(x_0) \eta(x_1) -  \eta(i x_1) \\
&= i T(x_0) -i  s(x_0) \eta(x_1) -  \eta(i x_1),
\end{aligned} $$ and thus $i  s(x_0) \eta(x_1) =  -  \eta(i x_1)$. Taking norms we arrive to $|s(x_0)|\  \| \eta(x_1) \| $ $=  \|-  \eta(i x_1)\|$ $= \|\eta(x_1)\|$ $\neq 0,$ and consequently $s(x_0) \in \{\pm 1\}$ and $\eta (x_0) =0$ (recall that $1 = \|x_0\|^2 = s(x_0)^2 + \|\eta(x_0)\|^2$), which is impossible. Therefore $x_1$ is a point of pure complex type for $T$.  \Cref{p the existence of a pure complex point implies that all are of pure complex type} proves that $T(i z) = i s(x_0) T(z)$ for all $z\in H$, and thus $T$ is complex-linear or conjugate-linear. 
\end{proof}

\begin{remark}\label{r: T is not complex-linear nor conjugate-linear} Under the hypotheses of \Cref{t characterization of complex-linear isometries OP}, it can be easily seen from the just quoted theorem and \Cref{p: mixed complex type point} that the following statements are equivalent: \begin{enumerate}[$(a)$]
		\item $T$ is not complex-linear or conjugate-linear.
		\item Every point in $H$ is of mixed complex type for $T$.
		\item There exists a non-zero point in $H$ which is of mixed complex type for $T$.
		\item There exists a non-zero point $x_0\in H$ such that $i T(x_0) \notin T(H)$.
	\end{enumerate}
\end{remark}

The next result is a consequence of \Cref{p real-linearity and distance preservation} and the previous theorem. 

\begin{theorem}\label{t: characterization of comple linear OP} Let $H$ and $K$ be complex inner product spaces such that $\dim (H)\ge 2$. Suppose that $A: H\to K$ is a non-zero additive mapping preserving Euclidean orthogonality. Then the following statements are equivalent:
\begin{enumerate}[$(a)$]
\item $A$ is complex-linear or conjugate-linear.
\item For every $z\in H$ we have $A(i z)  \in \{\pm i A(z)\}$.
\item There exists a non-zero point $z\in H$ such that $A(i z)  \in \{\pm i A(z)\}$.
\item There exists a non-zero point $z\in H$ such that $i A(z) \in A(H)$.
\end{enumerate}\smallskip

\noindent Furthermore, $A$ neither is complex-linear nor conjugate-linear if, and only if, there exists a non-zero $x\in H$ such that $i A(x)\notin A(H)$ (equivalently, for every non-zero $x\in H$, $i A(x)\notin A(H)$). 
\end{theorem}

Surprisingly, the algebraic characterization of those additive maps preserving orthogonality between complex inner product spaces which are complex-linear or conjugate-linear established in \Cref{t: characterization of comple linear OP} allows us to rediscover, via an alternative argument, the main result in \cite{LLP2024AOP}, namely, every additive mapping with dense image preserving orthogonality between complex inner product spaces must be complex-linear or conjugate-linear (cf. \cite[Theorem 2.1]{LLP2024AOP}).  

\begin{corollary}\label{c: if A has norm-dense range every point in the domain must be of mixed type} Let $H$ and $K$ be complex inner product spaces such that $\dim (H)\ge 2$. Suppose that $A: H\to K$ is a non-zero additive mapping preserving Euclidean orthogonality. Then if $A$ has norm-dense range, every point in $H$ is of pure complex type for the real-linear isometry associated with $A$ via \Cref{p real-linearity and distance preservation}, and $A$ must be a positive scalar multiple of a complex-linear or a conjugate-linear isometry.  
\end{corollary}   

\begin{proof} Suppose that $A$ has norm-dense range. Let us find, via \Cref{p real-linearity and distance preservation}, a positive real $\gamma$ such that $T= \gamma^{-1} A$ is a real-linear isometry preserving orthogonality. It follows from the hypotheses that $T$ has dense range. Clearly $T$ is complex-linear or conjugate-linear if, and only if, $A$ is. If $A$ (equivalently $T$) is not complex-linear or conjugate-linear every point in $H$ is of mixed complex type for $T$. Let us fix a norm-one vector $x_0\in H$. Note that $\eta(x_0)\neq 0$ in $K$.\smallskip
	
For any other $x\in H$ we can find $\alpha,\beta \in \mathbb{C}$ and a norm-one $x_1\in H$ such that $\langle x_0 | x_1\rangle =0,$ $|\alpha|^2 + |\beta|^2 =\|x\|^2$, and $x = \alpha x_0 + \beta x_1$. The sets $\{T(x_0), \eta(x_0)\}$ and  $\{T(x_1), \eta(x_1)\}$ are orthogonal in $K$ (cf. \Cref{p: mixed complex type point}). Consequently $T\left(\beta  x_1\right)\perp T(x_0), \eta(x_0), T(i x_0)$. The inequality $$\begin{aligned}
\| i \eta (x_0) - T(x) \|^2 &= \| i \eta (x_0) - \Re\hbox{e}(\alpha) T\left( x_0 \right) -\Im\hbox{m} (\alpha) T\left( i x_0 \right)- T\left(\beta  x_1\right)\|^2 \\
&= \| i \eta (x_0) - \Re\hbox{e}(\alpha) T\left( x_0 \right) -\Im\hbox{m} (\alpha) T\left( i x_0 \right)\|^2 +\| T\left(\beta  x_1\right)\|^2 \\
&\geq \| i \eta (x_0) - \Re\hbox{e}(\alpha) T\left( x_0 \right) -\Im\hbox{m} (\alpha) T\left( i x_0 \right)\|^2 \\
&= \| i \eta (x_0) - \Re\hbox{e}(\alpha) T\left( x_0 \right) -\Im\hbox{m} (\alpha) i s(x_0) T\left( x_0 \right) -\Im\hbox{m} (\alpha) \eta (x_0) \|^2 \\
& = |i -\Im\hbox{m} (\alpha)|^2 \|\eta (x_0)\|^2 + |\Re\hbox{e}(\alpha) -\Im\hbox{m} (\alpha) i s(x_0)|^2 \\
&= \|\eta (x_0)\|^2 + \Im\hbox{m} (\alpha)^2 \|\eta (x_0)\|^2 + \Re\hbox{e}(\alpha)^2 + \Im\hbox{m} (\alpha)^2 s(x_0)^2\\
&= \|\eta (x_0)\|^2 + \Im\hbox{m} (\alpha)^2 + \Re\hbox{e}(\alpha)^2 \geq \|\eta (x_0)\|^2 >0,
	\end{aligned}$$ assures that $i \eta (x_0)$ is not in the norm-closure of $T(H)$, which contradicts that $T$ has dense range. Therefore, every point in $H$ is of pure complex type for $T$, and \Cref{t characterization of complex-linear isometries OP} proves that $T$ is complex-linear or conjugate-linear.
\end{proof}

Having norm-dense range is a sufficient condition to guarantee that an additive mapping preserving orthogonality between complex inner product spaces is complex-linear or conjugate-linear. However, as the next corollary shows, it is not necessary.

\begin{corollary}\label{c automatic complex-linearity or conjugate linearity for finite ranges} Let $H$ and $K$ be finite dimensional complex Hilbert spaces with dim$(K)< 2\hbox{dim}(H)$ and dim$(H)\geq 2$. Then every additive mapping $A: H\to K$ preserving orthogonality is complex-linear or conjugate-linear.  
\end{corollary}   

\begin{proof} Arguing by contradiction, we assume that $0\neq A$ is not complex-linear or conjugate-linear.  
Let $\{e_1,\ldots, e_m\}$ be an orthonormal basis of $H$. A combination of \Cref{p real-linearity and distance preservation}, \Cref{r: T is not complex-linear nor conjugate-linear} and \Cref{p: mixed complex type point} shows that the vectors in the set $\{A(e_1),\ldots, A(e_m), \eta(e_1),\ldots, \eta(e_m) \}$ are mutually orthogonal in $K$, which is impossible since dim$(K)< 2 m$.
\end{proof}


We complement the conclusions above by throwing new light into the study of orthogonality preservers which are not complex-linear or conjugate-linear for this purpose we shall work with complex Hilbert spaces. A first example can be given as follows: Let $\tilde{H}$ and $\tilde{K}$ be complex Hilbert spaces for which we can find orthonormal systems $\{e_j : j\in \Gamma\}\subset H$ and $\{k_j, \tilde{k}_j: j\in \Gamma\}\subset K$, where the first one is a basis of $\tilde{H}$. Let us fix non-zero real numbers $s,c\in (0,1)$ satisfying $c^2 + s^2 =1$, and a selection $\sigma: \Gamma\to \{\pm 1\}$. We define a real-linear mapping $T: \tilde{H}\to \tilde{K}$ given by $T(e_j) = k_j$ and $T(i e_j) = i s k_j + \sigma(j) c \tilde{k}_j$  ($j\in \Gamma$), that is, \begin{equation}\label{eq general form of a real-linear isometry op not c linear} T\left(\sum_{j\in \Gamma} \alpha_j e_j \right) =  \sum_{j\in \Gamma} \Re\hbox{e}(\alpha_j) k_j +i  s \Im\hbox{m}(\alpha_j) {k}_j + \sigma(j) c \Im\hbox{m}(\alpha_j) \tilde{k}_j.
\end{equation} It is easy to check that $T$ is a real-linear isometry preserving orthogonality which is not complex-linear or conjugate-linear. \smallskip

The example exhibited above is essentially the unique possible choice in the case of complex Hilbert spaces. Namely, suppose $T: \tilde{H}\to \tilde{K}$ is a real-linear isometry preserving orthogonality between two complex Hilbert spaces which is not complex-linear or conjugate-linear. Let $\{e_j : j\in \Gamma\}$ be an orthonormal basis of $\tilde{H}$. By the hypotheses on $T$ and \Cref{p: mixed complex type point}, the set $\{T(e_j), \frac{\eta(e_j)}{\|\eta(e_j)\|}  : j\in \Gamma\}$ is an orthonormal system in $\tilde{K}$, and since $\{e_j : j\in \Gamma\}$ is a basis of $\tilde{H}$, $T$ admits a representation of the form in \eqref{eq general form of a real-linear isometry op not c linear} with $k_j = T(e_j),$ $\tilde{k_j}=\frac{\eta(e_j)}{\|\eta(e_j)\|}$, $\sigma(j) = 1$ for all $j$, $s= \langle T(i e_j)| i T(e_j)\rangle$, and $c = \|\eta(e_j)\|$ (which do not depend on $j$ by \Cref{p: mixed complex type point}). \smallskip

We can now define a real-linear mapping $R : \tilde{K}\to \tilde{K}$ given by $R (T(e_j)) = T(e_j),$ $R (i T(e_j)) = i s T(e_j),$ and $R \left(\frac{\eta(e_j)}{\|\eta(e_j)\|}\right) = i c T(e_j).$ Clearly $R$ is not an isometry, and does not preserve orthogonality ($\eta(e_j)\perp T(e_j)$ but $R(\eta(e_j))\not\perp R(T(e_j))$). However, it is not hard to check that the composed mapping $RT: \tilde{H}\to \tilde{K}$ is a complex-linear isometry preserving orthogonality.\smallskip

We have proved the following:

\begin{lemma}\label{l: existence of a liner transformation making T complex-linear} Let $T: \tilde{H}\to \tilde{K}$ be a real-linear isometry preserving orthogonality between two complex Hilbert spaces. Then there exists a bounded real-linear mapping $R: \tilde{K}\to \tilde{K}$ satisfying that $RT: \tilde{H}\to \tilde{K}$ is a complex-linear isometry preserving orthogonality. 
\end{lemma}

\begin{proof} If $T$ is not complex-linear or conjugate-linear we apply the previous arguments, while if $T$ is conjugate-linear, we compose with the conjugation on $\tilde{K}$ obtained by scalar multiplication by $i$.  
\end{proof}

If $T: \tilde{H}\to \tilde{K}$ is a real-linear isometry between two complex Hilbert spaces, we can regard $H$ and $K$ as real-Hilbert spaces with inner product $(x|y)= \Re\hbox{e}\langle x| y\rangle$ denoted by $\tilde{H}_{\mathbb{R}}$ and $\tilde{K}_{\mathbb{R}}$. If $\{e_j : j\in \Gamma\}$ be an orthonormal basis of $\tilde{H}$ is an orthonormal basis for $\tilde{H}$, the set $\{e_j, i e_j : j\in \Gamma\}$ is an orthonormal basis of $\tilde{H}_{\mathbb{R}}$. A similar construction to that given above proves the existence of a real-linear isometry $R: \tilde{K}_{\mathbb{R}} \to \tilde{K}_{\mathbb{R}}$ such that $R T$ is a complex-linear isometry preserving orthogonality.\smallskip

In order to consider orthogonality preservers between inner product spaces we need the following lemma. 

\begin{lemma}\label{l op lift to the completion} Let $A: H\to K$ be an additive orthogonality preserving mapping between two $($real or complex$)$ inner product spaces with dim$(H)\geq 2$. Then $A$ is real-linear and there exists a unique extension of $A$ to a real-linear orthogonality preserving map between the completions of $H$ and $K$.  
\end{lemma} 

\begin{proof} We can clearly assume that $A\neq 0$. \Cref{p real-linearity and distance preservation} implies the existence of a positive real number  $\gamma$, and a real-linear isometry preserving orthogonality $T: H\to K$ such that $A = \gamma T$. Let $\tilde{H}$ and $\tilde{K}$ denote the completions of ${H}$ and ${K},$ respectively. We can obviously extend (uniquely) $A$ and $T$ to real-linear maps $\tilde{A},\tilde{T}: \tilde{H} \to \tilde{K}$, where $\tilde{T}$ is a real-linear isometry and $\tilde{A} = \gamma \tilde{T}$. It remains to prove that $\tilde{A}$ and $\tilde{T}$ preserve orthogonality. The conclusion is clear when $H$ and $K$ are real inner product spaces. We therefore consider the case of complex inner product spaces. Let us take $z,w\in \tilde{H}\setminus\{0\}$ with $\langle z| w\rangle = 0$, and find two sequences $(z_n)_n, (w_n)_n$ in $H$ converging in norm to $z$ and $w$, respectively.  We can also assume that $z_n\neq 0$ for all $n$. Define $\hat{w}_n = w_n - \langle w_n | z_n\rangle \frac{z_n}{\| z_n\|^2}$. It is easy to see that $\langle \hat{w}_n | z_n\rangle =0$ for all $n$, and $(\hat{w}_n)\to w$.  By hypotheses, $\langle A(\hat{w}_n) | A(z_n)\rangle =0$ for all $n$, and by the continuity of $\tilde{A}$ we get $\langle \tilde{A}({w}) | \tilde{A}(z)\rangle =0$, which concludes the proof.
\end{proof}

We can actually prove a more general result.

\begin{lemma}\label{l op lift to the completion normed with dual strictly convex} Let $A: X\to Y$ be an additive mapping between two complex normed spaces preserving Birkhoff orthogonality. Suppose that dim$(X)\geq 2,$ $X$ admits a conjugation $($i.e. a conjugate-linear isometry of period-$2)$, and $X^*$ is a strictly convex space. Then $A$ is real-linear and there exists a unique extension of $A$ to a real-linear orthogonality preserving map between the completions of $X$ and $Y$.  
\end{lemma} 

\begin{proof} Proposition 2.2 in \cite{LLP2024AOP} implies that $A$ is a bounded real-linear mapping. We can clearly extend $A$ to a bounded real-linear mapping $\tilde{A}$ between the completions $\tilde{X}$ and $\tilde{Y}$ of $X$ and $Y$, respectively. It only remains to prove that $\tilde{A}$ preserves Birkhoff orthogonality. \smallskip
	
Since $X^*$ and $\tilde{X}^*$ are isometrically isomorphic, it follows from the assumptions on ${X}^*$ and a classic result in functional analysis (see \cite[Theorem in page 172]{HolmesBook75}) that $\tilde{X}$ is smooth. The latter implies that every bounded sequence in $\tilde{X}^*$ has a weak$^*$ converging subsequence (\cite[Theorem 1]{HagSul1980}).\smallskip

Let us find $x,y \in \tilde{X}$	with $x\perp_{B} y$, which is equivalent to the existence of a norm-one functional $\phi_0\in \tilde{X}^*$ satisfying $\phi_0 (x) =\|x\|$ and $\phi_0 (y)  =0$. Let us pick sequences $(x_n)_n$ and $(y_n)_n$ in $X$ converging in norm to $x$ and $y$, respectively. We can assume that $x,y,x_n,y_n\neq 0,$ for all $n\in \mathbb{N}$. For each $n\in \mathbb{N}$, find, via Hahn-Banach theorem, norm-one functional $\phi_n\in \tilde{X}^*$ satisfying $\phi_n (x_n) = \|x_n\|$, for all natural $n$. The element $z_n = y_n-\frac{\phi_n (y_n)}{\|x_n\|} x_n$ lies in $X$ with $x_n\perp_{B} z_n$ ($n\in \mathbb{N}$). We know from the previous paragraph that there is a subsequence $(\phi_{n_k})_k$ converging to some $\phi_1$ in the weak$^*$ topology of $\tilde{X}^*$. The identity $$\|x_{n_k}\| = \phi_{n_k} (x_{n_k}) = \phi_{n_k} (x)-\phi_{n_k} (x-x_{n_k}),$$ and the properties of the sequences $(\phi_{n_k})_k$ and $(x_{n_k})_k$ assure that $\phi_1 (x) = \|x\|$, and thus $\|\phi_1 \| =1$. It follows from the smoothness of $\tilde{X}$ that $\phi_1 =\phi_0$. Therefore $$  (\phi_{n_k} (y_{n_k}))_k = (\phi_{n_k} (y)-\phi_{n_k} (y-y_{n_k}))_k\to \phi_0 (y) =0,$$ which consequently proves that $(z_{n_k})_k\to y$ in norm.\smallskip 

Finally, since $A$ preserves Birkhoff orthogonality, $A(x_{n_k}) \perp_{B} A(z_{n_k})$ for all $k\in \mathbb{N}$, that is, $$\| A(x_{n_k}) + \lambda A(z_{n_k}) \|  \geq \|A(x_{n_k})\|, \hbox{ for all } \lambda\in \mathbb{C}.$$ Taking limits in $k$ and considering the continuity of $\tilde{A}$, we deduce that $\left\| \tilde{A} (x) + \lambda \tilde{A}(y) \right\|  \geq \left\|\tilde{A}(x)\right\|,$ for all $\lambda\in \mathbb{C}.$
\end{proof}

Suppose $A: X\to Y$ is an additive mapping between two real normed spaces preserving Birkhoff orthogonality, with dim$(X)\geq 2.$ Then W\'{o}jcik's theorem \cite[Theorem 3.1]{Wojcik2019} proves that $A$ is real-linear, and thus there exists a unique extension of $A$ to a real-linear orthogonality preserving map between the completions of $X$ and $Y$, since in the real setting every positive scalar multiple of a real-linear isometry preserves Birkhoff orthogonality.\smallskip

We conclude this note with the following consequence of the previous two Lemmas~\ref{l: existence of a liner transformation making T complex-linear} and \ref{l op lift to the completion}.

\begin{corollary}\label{c: final} Let $A: H\to K$ be an additive orthogonality preserving mapping between two complex inner product spaces with dim$(H)\geq 2$, and let $\tilde{K}$ denote the completion of ${K}$. Then there exists a bounded real-linear mapping $R: \tilde{K}\to \tilde{K}$ satisfying that $R T: {H}\to \tilde{K}$ is a complex-linear positive scalar multiple of a complex-linear isometry preserving orthogonality. 
\end{corollary}

\smallskip

\textbf{Acknowledgements} L. Li was supported by National Natural Science Foundation of China (Grant No. 12171251). A.M. Peralta supported by grant PID2021-122126NB-C31 funded by  MICIU/AEI/ 10.13039/501100011033 and by ERDF/EU, by Junta de Andalucía grant FQM375, \linebreak IMAG--Mar{\'i}a de Maeztu grant CEX2020-001105-M/AEI/10.13039/501100011033 and (MOST) Ministry of Science and Technology of China grant G2023125007L. \smallskip


\subsection*{Statements and Declarations} 

All authors declare that they have no conflicts of interest to disclose.

\subsection*{Data availability}

There is no data associate for this submission.

\end{document}